\documentclass{article}

\newtheorem{thm}{Theorem}[section]

\newtheorem{cor}[thm]{Corollary}
\newtheorem{ex}[thm]{Example}

\newtheorem{remark}[thm]{Remark}

\newcommand\ack{\section*{Acknowledgement.}}

\newcommand{\A}{\mbox{$\cal A$}}
\newcommand{\F}{\mbox{$\cal F$}}
\newcommand{\C}{\mbox{$\cal C$}}
\newcommand{\G}{\mbox{$\cal G$}}
\newcommand{\bbI}{{\rm I\hspace{-0.8mm}I}}
\newcommand{\bbR}{{\rm I\hspace{-0.8mm}R}}

\newcommand{\bbN}{{\rm I\hspace{-0.8mm}N}}
\newcommand{\bbP}{{\rm I\hspace{-0.8mm}P}}
\newcommand{\bbE}{{\rm I\hspace{-0.8mm}E}}
\newcommand{\matR}{{\bbR}}
\newcommand{\matP}{{\bbP}}
\newcommand{\matI}{{\bbI}}
\newcommand{\matE}{{\bbE}}
\newcommand{\koniecmat}{\,}

\newcommand{\ogon}[1]{\overline{#1}}
\newcommand{\etal}{{\it et al. }}

\begin{document}

\title{Poisson limits for empirical point processes} 

\author{Andr\'{e} Dabrowski*\\
Department of Mathematics and Statistics\\
Ottawa University\\
585 King Edward Av.\\
K1N 6N5 Ottawa, Canada\\
ardsg@science.uottawa.ca
 \and Gail Ivanoff\\
Department of Mathematics and Statistics\\
Ottawa University\\
585 King Edward Av.\\
K1N 6N5 Ottawa, Canada\\
givanoff@science.uottawa.ca
 \and Rafa{\l} Kulik\\
School of Mathematics and Statistics\\
F07, University of Sydney\\
NSW 2006, Australia
 }

\maketitle



\begin{abstract}
Define the scaled empirical point process on an independent and
identically distributed sequence $\{Y_i:\ i\le n\}$ as the random
point measure with masses at $a_n^{-1} Y_i$. For suitable $a_n$ we
obtain the weak limit of these point processes through a novel use
of a dimension-free method based on the convergence of
compensators of multiparameter martingales. The method extends
previous results in several directions. We obtain limits at points
where the density of $Y_i$ may be zero, but has regular variation.
The joint limit of the empirical process evaluated at distinct
points is given by independent Poisson processes. These results
also hold for multivariate $Y_i$ with little additional effort.
Applications are provided both to nearest-neighbour density
estimation in high dimensions, and to the asymptotic behaviour of
multivariate extremes such as those arising from bivariate normal
copulas.
\end{abstract}

\noindent{\bf Keywords:} multiparameter martingales; point
processes; density estimation; multivariate extremes; local
empirical processes

\noindent{\bf Running title:} Empirical point processes
\newpage
\section{Introduction}

Point processes and their limits arise naturally in many areas of
statistics, and have a number of applications ranging from
survival analysis to spatial statistics. Point processes also
arise in probability theory as limits for extreme value processes,
in studying limits of sums of stable non-Gaussian variables and in
queuing models. Of course the Poisson process is a fundamental
concept in martingale theory. Weak convergence of the empirical
point process underlies many applications, and this paper employs
the relatively recent area of multiparameter martingales to
establish a novel and unified approach to proving such limits for
scaled empirical point processes. Although various elegant and
powerful methods have been developed for particular
 classes of problems, the generalized martingale approach provides an extremely  simple,
dimension-free method of addressing a variety of old and new
distributional questions.

Given a random sample of random vectors $\{Y_i:\ i\le n\}$ in
${\matR}^d$ and a suitable class of sets $\{ A\}$, the empirical
point process is defined by
$$N^{(n)}_A=\sum_{i=1}^n\bbI_{\{Y_i\in A\}}.$$
As noted above, the weak convergence of $N^{(n)}_A$ has been
extensively studied using a variety of methods. In particular, a
strong approximation approach can be used to establish weak
convergence of the local empirical process (see Einmahl, 1997, and
the references therein):
$$
L_{n,x}(t)=\frac{1}{n}\sum_{i=1}^n \bbI_{\{Y_i\in [x-ta_n,x+ta_n]\}}
,\ \ \ t\in [0,1],
$$
where now the $Y_i$'s are univariate. If the sequence of
constants, $a_n$, is appropriately chosen then the limit process
is homogeneous Poisson. However, this strong Poisson approximation
is difficult to implement (or at least cannot be extended
directly)
if one wants to study the joint behaviour of
$$
(L_{n,x_1}(\cdot),\ldots,L_{n,x_m}(\cdot)),
$$
 i.e. when estimating the
density of $Y_1$ simultaneously at $(x_1,\ldots,x_m)$ (see Section
4.1). Even in the Gaussian case, where simultaneous approximation by
independent Wiener processes is known, Deheuvels \etal (2000) points
out that a major technical difficulty arises in proving independence
at separate $x_i$.

 The aim of this paper is to develop a general and natural
approach to \textit{weak} Poisson limits for empirical point
processes. It is based on the multiparameter martingale theory of
Ivanoff and Merzbach (2000) and requires only the simple computation
of so-called
*-compensators to identify Poisson limits for scaled empirical point
processes. The compensator method exploited here is particularly
attractive in that it is independent of the dimension of the
underlying random vectors, and so easily generalizes results from
the univariate to the multivariate case. In addition, the martingale
approach allows one to handle  the joint behaviour at multiple
points  with ease through a judicious definition of the associated
history (filtration). In particular, we shall show that the
asymptotic behaviour of the local empirical process at distinct
points $x_1,\ldots,x_m$ can be described by independent Poisson
processes, an intuitive but otherwise technically challenging
result.

The method has additional benefits. First, only (multivariate)
regular variation of the density $f$ of $Y_1$ is required, and the
limits are explicitly written in terms of $f$. Indeed, we can
discover the appropriate scaling constants even when $f$ is
regularly varying but $f(0)=0$, i.e. a case with inhomogeneous
Poisson limits excluded in Borisov (2000).
characterize the distributional behaviour of joint extremes for
different bivariate copulas. This recovers Einmahl (1997,
Corollaries 2.4 and 2.5) where the limit Poisson process has a
product mean measure, but also extends to more complex cases
(Corollary~\ref{corcopula}). In particular, we can identify extreme
value limits for copulas with asymptotically dependent multivariate
extremes more simply than methods employing multivariate regular
variation (c.f. Resnick, 1987).

The paper is structured as follows. The next section will review key
elements of the theory of multi-parameter martingales and in
particular, the use of
*-compensators in proving weak convergence of a sequence of point
processes. Section~\ref{sec.quantiles} defines the scaled empirical
point process generated by a sample and establishes point process
limits for such processes. This proceeds in steps from the classical
non-negative and univariate case (yielding limits similar to those
 for extreme value processes), to the multivariate and
multidimensional cases. In each case the proof simplifies to the
straightforward calculation of *-compensators, and highlights the
universality of the martingale approach. Section~\ref{sec.applic}
on applications illustrates the utility of our results by
establishing for the first time weak limits for nearest-neighbour
estimates of joint densities (again at several points
simultaneously), and by providing new extreme value limits for
multivariate copulas.

\section{Notation and background: Point processes and martingale
methods}

We provide a brief introduction to point processes and martingale
methods indexed by general Euclidean spaces using the set-indexed
framework introduced in Ivanoff and Merzbach (2000). We need
definitions mimicking those for martingales indexed by ${\matR}_+$.

Set $T={\bbR}^d$ or ${\bbR}_+^d$, and ${\cal A}=\{A_t=[0,t]:t\in
T\}\cup\{\emptyset\}$, where we interpret $[0,t]$ in the obvious way
if $t\not\in {\bbR}_+^d$.  Set-inclusion on ${\cal A}$ induces a
partial order, $\preceq$, on $T$: $s\preceq t$ if and only if
$A_s\subseteq A_t$. This is not the usual partial order on $\bbR^d$:
e.g. $\{0\}$ is the (unique) minimal element, and all quadrants are
equipped with their own partial order. In particular, if
$T={\matR}$, points with different signs are incomparable. This
special structure permits us to define a $2^d$-sided martingale
theory.

The semi-algebra\index{semi-algebra} ${\cal C}$ \index{${\cal C}$}
is the class of all subsets of $T$ of the form
$$C=A\setminus B,~A\in {\cal A},~ B\in {\cal A} (u),$$ where ${\cal A}
(u)$ denotes the class of sets which are finite   unions of sets
from ${\cal A}$. Let $(\Omega,\F , P)$ be any complete probability
space. A filtration indexed by ${\cal A}(u)$ is a class $ \{ {
\cal F } _A :A\in{\cal A}(u) \} $ \index{filtration}\index{${\cal
F}_A$} of complete sub-$\sigma$-fields of $\F$ where $ \forall A ,
B \in { \cal A }(u) $, $ { \cal F } _ A \subseteq { \cal F } _ B $
if $ A \subseteq B $, and (Monotone outer-continuity) $ { \cal F }
_{ \cap A _ i } = \cap { \cal F } _ {A_i}$ for any decreasing
sequence $(A_i)$ in   $ { \cal A }(u)$ such that $\cap_iA_i\in
{\cal A}(u)$. For consistency, we define $\F_T=\F$. We may
associate $\sigma$-algebras with sets in ${\cal C}$: for $C\in
\C\setminus {\cal A}$, let
$\G _C^*=\vee_{B\in {\cal A}(u), B\cap C=\emptyset}\F_B$, and for
$A\in {\cal A}, A\neq \emptyset$, define $\G_A^*=\F_{\emptyset}$. A
\emph{(${\cal A}$ -indexed) stochastic process} $X=\{X_A:A\in{\cal
A}\}$  is a collection of random variables indexed by ${\cal A}$,
and is \emph{adapted} if $X_A$ is $\F_A$-measurable for every $A\in
{\cal A}$. By convention, $X_{\{0\}}=0$.

A process $X:{\cal A}\rightarrow \bbR  $ is \emph{increasing} if for
every $\omega\in \Omega$, the function $X_{\cdot}(\omega)$ can be
extended to a finitely additive function on $\C$ satisfying
$X_{\{0\}}(\omega)=0$ and $X_C(\omega)\geq    0,~\forall C\in \C$,
and such that if $(A_n)$ is a decreasing sequence of sets in $\A
(u)$ such that $\cap_nA_n\in\A(u)$, then $\lim_n  X_{A_n}(\omega
)=X_{\cap_nA_n}(\omega)$.  A process $N = \{N_A , A \in\cal{A}\}$ is
a \emph{point process} if it is an increasing process taking its
values in ${\bbN}$, and almost surely for any $t\in T$,
$N_{\{t\}}=0$ or 1.
Note that if $N$ is a point process on $T=\bbR  $, then
$N_t:=N_{[0,t]}$ (for $t$ positive or negative) and not $N_{(-\infty
,t]}$.  As expected, $N$ is a Poisson process on $T$ with mean
measure $\Lambda$ if $N$ is a point process where $N_C\sim
\mbox{Poisson}, \Lambda_C, ~\forall C\in \C$, and whenever
$C_1,...,C_n\in \C$ are disjoint, $N_{C_1},...,N_{C_n}$ are
independent. If $\Lambda$ is absolutely continuous with respect to
Lebesgue measure, its density $\lambda$ is called the intensity of
the Poisson process.

An integrable process $M = \{M_A , A \in\cal{A} \}$ is called a {\it
pseudo-strong martingale} if for any $C \in\cal{C}$, $ E[ M_C |
{\cal G}_{C}^{*} ]$ = $0$.
The process $\overline{X}$ is a *-{\it compensator} of $X$ if it is
increasing and the difference $X - \overline{X} $ is a pseudo-strong
martingale.
The asymptotic behaviour of a sequence of point processes may be
determined by *-compensators as shown in the following theorem
specializing Theorem 8.2.2 and Corollary 8.2.3 of Ivanoff and
Merzbach (2000) to multivariate point processes on $T={\matR}  ^d$
or ${\matR}_+^d$.

To state this theorem, we consider $k$ point processes
$N(1),...,N(k)$ all adapted to a common $\A (u)$-indexed filtration
$\{\F_A\}$ and so that with probability one, none of the processes
have a jump point in common. The $k$-variate point process
$\stackrel{\rightarrow}{N}$ is defined by
$\stackrel{\rightarrow}{N}_A=\langle N_A(1),...,N_A(k)\rangle$ and
has ($k$-variate) *-compensator
$\stackrel{\rightarrow}{\Lambda}=(\Lambda(1),...,\Lambda (k))$ if
$\Lambda(i)$ is a *-compensator for $N(i)$ with respect to the
common filtration $\{\F_A\}$.

In what follows,``$\longrightarrow_P$" denotes convergence in
probability and ``${\longrightarrow_{\cal D}}$" denotes convergence
in both finite dimensional distribution and in distribution in the
Skorokhod topology if $T=\bbR  ^d_+$ (identifying $N^{(n)}_t
$(respectively, $N_t$) with $N^{(n)}_{A_t}$ (respectively, $N_{A_t}
$)). We remark that the Skorokhod topology may be extended to all of
the quadrants in $\bbR  ^d$ on the space of ``outer-continuous
functions with inner limits", and the convergence in the theorem
above holds in this case as well. In the sequel, convergence in the
Skorokhod topology will be interpreted in this way.

\begin{thm}\label{multichannel} Let
$(\stackrel{\rightarrow}{N}^{(n)})$ be a sequence of $k$-variate
point processes on $T$ adapted to a filtration $\{{\cal F}_A\}$
and $(\stackrel{\rightarrow}{\Lambda} ^{(n)})$ a sequence of
corresponding *-compensators. Suppose that for each $A\in \A$ and
$i=1,...,k$ the sequences $(N_A^{(n)}(i))$ and
$(\Lambda_A^{(n)}(i))$ are uniformly integrable and that
$\Lambda^{(n)}_A(i)\longrightarrow_P \Lambda_A(i)$ where
$\Lambda(i)$ is a deterministic   measure on $T$ absolutely
continuous  with respect to Lebesgue measure. Then
$\stackrel{\rightarrow}{N}^{(n)}{\longrightarrow_{\cal
D}}\stackrel{\rightarrow}{N}$, where
$\stackrel{\rightarrow}{N}=\langle N(1),...,N(k)\rangle$ and
$N(1),...,N(k)$ are independent  Poisson processes with mean
measures $ \Lambda (1),...,\Lambda (k)$, respectively.

\end{thm}

\noindent {\bf Proof}: The proof of this theorem is a
straightforward generalization of the techniques used in the proof
of Theorem 8.2.2 in Ivanoff and Merzbach (2000) along with an
application of Watanabe's characterization of the $k$-variate
Poisson process on $\bbR _+$, see Br\'{e}maud (1981, Theorem
T6).\hfill

We conclude this section by defining empirical point processes on
$T$ and stating their *-compensators. Let $Y$ be a $T$-valued random
variable with continuous distribution function $F$. The {\em single
jump} point process $J=\{J_A=\bbI_{\{Y\in A\}}:A\in \A\}$ has
*-compensator
$$\overline{J}_A=\int_A\bbI_{\{u\in
A_Y\}}(F(E_u))^{-1}dF(u)$$
  with respect to its minimal filtration, where $E_t=\{t'\in T: t\preceq
t' \}$ (cf. Ivanoff and Merzbach (2000)).  Now, suppose that
$Y_1,..., Y_n$ are i.i.d. with distribution $F(t)={\matP}(Y_i\le t)$
and let $\F=\vee_{i=1}^n \F^{(i)}$ where $\F^{(i)}$ is the minimal
filtration generated by the single jump process associated with
$Y_i$. Then the \emph{empirical point process} $N^{(n)}$ defined by
$$N^{(n)}_A=\sum_{i=1}^n\bbI_{\{Y_i\in A\}}$$
has *-compensator $\Lambda^{(n)}$ where
\begin{equation}\label{*comp}\Lambda^{(n)}_A=\sum_{i=1}^n\int_A\bbI_{\{u\in
A_{Y_i}\}}(F(E_u))^{-1}dF(u)\koniecmat .
\end{equation}
\begin{ex}{\rm If $T={\matR}_+$ then (\ref{*comp}) reads as
follows: $A=A_t=[0,t]$ for some $t\ge 0$, $\preceq$ is just $\le$,
the standard ordering, $E_u=[u,\infty)$. By noting that
$F(E_u)={\matP}(Y_i\ge u)=:\ogon{F}(u)$ we have
\begin{equation}\label{ex.1dim.positive}
\Lambda_t^{(n)}:=\Lambda_{A_t}^{(n)}=
\int_{A_t}\sum_{i=1}^n{\matI}\{Y_i\in E_u\}\frac{dF(u)}{F(E_u)}=
\int_0^t\sum_{i=1}^n{\matI}\{Y_i\ge
u\}\frac{dF(u)}{\ogon{F}(u)}\koniecmat\; .
\end{equation}
 }
\end{ex}
\begin{ex}\label{exam.1dim.general}{\rm If $T={\matR}$ then $A_t=[0,t]$ or $A_t=[t,0]$
depending on the sign of $t$. We have $s\preceq t$ if $0\le s\le t$,
or $t\le s\le 0$, where $\le$ is the standard order on ${\matR}$.
Points with different signs are incomparable. The sets $E_u$ will be
either $[u,\infty)$ or $(-\infty,u]$ depending the sign of $u$. If
$u>0$ then as above $F(E_u)=\ogon{F}(u)$, otherwise, if $u<0$, then
$F(E_u)=F(u)$. Now the *-compensator is given by
(\ref{ex.1dim.positive}) if $t\ge 0$ and if $t<0$,
\begin{equation}\label{ex.1dim.general}
\Lambda_t^{(n)}:=\Lambda_{A_t}^{(n)}=
\int_{A_t}\sum_{i=1}^n{\matI}\{Y_i\in E_u\}\frac{dF(u)}{F(E_u)}=
\int_t^0\sum_{i=1}^n{\matI}\{Y_i\le u\}\frac{dF(u)}{F(u)}\koniecmat
.
\end{equation}
}
\end{ex}
\begin{ex}{\rm Let $Y_i=(Y_{i1},Y_{i2})$, $i=1,\ldots,n$. If $T={\matR}_+^2$ then
(\ref{*comp}) reads as follows: $A=A_{\bf t}=[0,t_1]\times [0,t_2]$
for some ${\bf t}=(t_1,t_2)$, $E_{\bf u}=\{{\bf t'}:t_i'\ge
u_i,i=1,2\}$, ${\bf u}=(u_1,u_2)$. By noting that $F(E_{\bf
u})={\matP}(Y_{i1}\ge u_1,Y_{i2}\ge u_2)=:\ogon{F}({\bf u})$ we have
$$
\Lambda_{\bf t}^{(n)}:=\Lambda_{A_{\bf t}}^{(n)}= \int_{A_{\bf
t}}\sum_{i=1}^n{\matI}\{Y_{i1}\ge u_1,Y_{i2}\ge u_2\}\frac{dF({\bf
u})}{\ogon{F}({\bf u})}\koniecmat\; .
$$
To extend this example to ${\matR}^2$ we proceed as in Example
\ref{exam.1dim.general}, treating each quadrant separately.
 }
\end{ex}

\section{Poisson limits at quantiles}\label{sec.quantiles}
\subsection{Univariate case}\label{sec.univariate} We can use the
previous section to determine the limiting behaviour of empirical
point processes at quantiles. Consider a sequence $\{Y_n\}$ of
i.i.d. real-valued positive random variables with distribution
$F$. Assume now that $F(0)=0$ and that $F$ is regularly varying at
0 with index $\alpha>0$, i.e. for all $t\ge 0$,
\begin{equation}\label{def.univ.rv}
\lim_{x\searrow 0}\frac{F(xt)}{F(x)}=t^{\alpha} \koniecmat ,
\end{equation}
see e.g. Resnick (1987). This implies that for $x$ in a
neighbourhood of $0$, $F(x)=\ell(x)x^\alpha$. Here and in the sequel
$\ell$ is a slowly varying function at $0$ or at $\infty$ as
required, and it can be different at each appearance.

Let $a_n$ be such that $F(a_n)=n^{-1}$. This ensures that $a_n\sim
n^{-1/\alpha}\ell(n)$ for some function $\ell$ slowly varying at
$\infty$. Henceforth, we write $c_n\sim d_n$ if $\lim_{n\to\infty}
c_n/d_n=1$.

Since $a_n\to 0$ we have
\begin{equation}
nF(a_nt)  =  \frac{F(a_nt)}{F(a_n)}  \to
t^{\alpha}\label{reg.var.comp}\; .
\end{equation}
Define
$$N^{(n)}=\sum_{i=1}^n\delta_{a_n^{-1}Y_i}\; .$$
We have by (\ref{*comp})
\begin{equation}\label{comp.empirical}
\Lambda_t^{(n)}= \int_0^t\sum_{i=1}^n{\matI}\{a_n^{-1}Y_i\ge
u\}\frac{dF(a_nu)}{\ogon{F}(a_nu)}\koniecmat\; .
\end{equation}

We first reprove the well-known result (see e.g. Resnick (1987,
Proposition 3.21) concerning Poisson limits for empirical point
processes. An elegant argument can be applied (see e.g. Borisov,
2001, and the references therein) where the law of $N^{(n)}$ is
approximated (in the total variation sense and for each $n$
separately) by $Poi(\nu_n)$, the Poisson random measures with
$\nu_n(A)=nE (1_{n^{1/k}X_i\in A})$, where the $X_i$'s are uniform
on an appropriately chosen ball. If $n$ is sufficiently large,
strong approximation methods yield the coupling of the empirical
point processes to a single Poisson random measure, and weak
convergence follows.  However here we illustrate the martingale
approach since the proof easily generalizes to the multivariate
context and to establishing simultaneous limits at interior
quantiles of $F$.

\begin{thm}\label{theo.1dim}
Assume that $F(0)=0$ and (\ref{def.univ.rv}) holds. Then the
sequence $(N^{(n)})$ converges in distribution to $N$ in the
Skorokhod topology on $D[0,\infty)$ where $N$ is a Poisson process
with mean measure $\Lambda_t=t^{\alpha}$ (intensity
$\lambda(t)=\alpha t^{\alpha-1}$).
\end{thm}
{\it Proof.} Since $N^{(n)}$ is square integrable with bounded
second moments (uniformly in $n$), the conditions of Theorem
\ref{multichannel} will be satisfied if it is shown that the
sequence $(\Lambda_t^{(n)})$ given by (\ref{comp.empirical})
converges in $L_2$ to $t^{\alpha}$.
\begin{eqnarray*}
{\matE}[\Lambda_t^{(n)}] & = &
{\matE}\left[\int_0^t\sum_{i=1}^n{\matI}\{a_n^{-1}Y_i\ge
u\}\frac{dF(a_nu)}{\ogon{F}(a_nu)}\right]\\
& = & \sum_{i=1}^n\int_0^t{\matP}(a_n^{-1}Y_i\ge
u)\frac{dF(a_nu)}{\ogon{F}(a_nu)}
 =  \sum_{i=1}^n\int_0^t dF(a_nu)
 =  nF(a_n t)
 \to  t^{\alpha} \koniecmat
\end{eqnarray*}
by applying (\ref{reg.var.comp}). Using the independence assumption,
\begin{eqnarray*}
{\matE}[(\Lambda_t^{(n)})^2] & = & \int_0^t\int_0^t
\sum_{i=1}^n{\matE}\left[{\matI}\left\{Y_i\ge a_nu,Y_i\ge
a_nv\right\}\right]
\frac{dF(a_nu)dF(a_nv)}{\ogon{F}(a_nu)\ogon{F}(a_nv)}\\
&& +2\int_0^t\int_0^t \sum_{i<j}^n{\matE}\left[{\matI}\left\{Y_i\ge
a_nu,Y_j\ge a_nv\right\}\right]
\frac{dF(a_nu)dF(a_nv)}{\ogon{F}(a_nu)\ogon{F}(a_nv)}\\
& = & \int_0^t\int_0^t \sum_{i=1}^n{\matP}\left(Y_i\ge a_n(u\vee
v)\right)
\frac{dF(a_nu)dF(a_nv)}{\ogon{F}(a_nu)\ogon{F}(a_nv)}\\
&& +n(n-1)\int_0^t\int_0^t
dF(a_nu)dF(a_nv)\\
& = & n\int_0^t\int_0^t \sum_{i=1}^n{\matP}\left(Y_i\ge a_n(u\vee
v)\right)
\frac{dF(a_nu)dF(a_nv)}{\ogon{F}(a_nu)\ogon{F}(a_nv)}\\
&& +n(n-1)(F(a_nt))^2\koniecmat .
\end{eqnarray*}
Now, the first term converges to 0, because
\begin{eqnarray*}
\lefteqn{n\int_0^t\int_0^t \sum_{i=1}^n{\matP}\left(Y_i\ge a_n(u\vee
v)\right) \frac{dF(a_nu)dF(a_nv)}{\ogon{F}(a_nu)\ogon{F}(a_nv)}
}\\
& \leq  & n\int_0^t\int_0^t \frac{dF(a_nu)dF(a_nv)}{\ogon{F}(a_nu)}
 \le  n\left[\ogon{F}(a_nt)\right]^{-1}\left[F(a_nt)\right]^{2},
\end{eqnarray*}
and $\ogon{F}(a_nt)\rightarrow 1$ and $nF^2(a_nt)\rightarrow 0$ as
$n\rightarrow \infty$. So, $(\Lambda^{(n)}_t)$ converges in $L_2$
and therefore in probability.

{\rm We may extend Theorem~\ref{theo.1dim} to the entire line. $P_F$
is the probability measure associated with a distribution $F$. We
say that $F$ is \emph{regularly varying on the right (left)} at $u$
with index $\alpha$ ($\beta$) if
\begin{equation}\label{def.univ.shift.rv.measure}\lim_{ x\searrow
0}\frac{P_F((u,u+xt])}{P_F((u,u+x])}=t^\alpha\;\; \left(\lim_{
x\searrow 0}\frac{P_F((u+xt,u])}{P_F((u-x,u])}=|t|^\beta \right)\; ,
\end{equation}
for all $t>0$ and $t<0$, respectively. Clearly, if $F$ has support
$(0,\infty)$ then for $u=0$ the above condition reduces to
(\ref{def.univ.rv}).

If $F$ fulfills (\ref{def.univ.shift.rv.measure}), we shall choose
$a_n$ and $b_n$ so that
\begin{equation}\label{shifted.threshold}
F(a_n+u)-F(u)=n^{-1}\;,\;\; F(u)-F(-b_n+u)=n^{-1}\; .
\end{equation}
Fix $q\in(0,1)$ and set $x_q=F^{-1}(q)$ and assume that $F$ fulfills
(\ref{def.univ.shift.rv.measure}) at $u=x_q$. Let
\begin{equation}\label{eq2}
 N^{(n)}(q)=\sum_{i=1}^n\delta_{a_n^{-1}[Y_i-x_q]}I[Y_i\ge x_q]
+\sum_{i=1}^n\delta_{b_n^{-1}[Y_i-x_q]}I[Y_i< x_q].
 \end{equation}
The argument in the preceding proof can now be repeated for $t>0$
and $t<0$ to prove that $\Lambda^{(n)}_t$ converges in $L_2$ to
$t^ {\alpha}$ if $t>0$ and to $|t|^\beta$ if $t<0$. Let $G$ be the
distribution of $Y_i-x_q$. Thus, $G(s)=F(s+x_q)$. To see that the
norming sequences $a_n$ and $b_n$ are chosen appropriately, using
the same calculation as in the proof of Theorem \ref{theo.1dim},
we have for $t>0$
\begin{eqnarray*}
{\matE}\left[\Lambda_{t}^{(n)}\right] & = & n\int_{0}^{t}dG(a_n u) =
n[F(a_nt+x_q )-F(x_q)] \\& =& \frac{F(a_nt+x_q
)-F(x_q)}{F(a_n+x_q)-F(x_q)}
 \to  t^{\alpha}
\end{eqnarray*}
by the first part of (\ref{def.univ.shift.rv.measure}). Moreover,
bearing in mind Example \ref{exam.1dim.general}, we have for $t<0$,
\begin{eqnarray*}
{\matE}\left[\Lambda_{t}^{(n)}\right] & = & n\int_t^0dG(b_n u)
 =  -n[F(b_nt+x_q )-F(x+q)]\\
& = & \frac{F(b_nt +x_q)-F(x_q)}{F(-b_n+x_q)-F(x_q)}
 \to  |t|^{\beta}
\end{eqnarray*}
by the second part of (\ref{def.univ.shift.rv.measure}).
Theorem~\ref{multichannel} leads to the following Corollary.

\begin{cor}\label{theo.q}
Assume (\ref{def.univ.shift.rv.measure}).
Then $N^{(n)}(q)\longrightarrow_{\cal D}N$, where $N$ is a Poisson
process on $\bbR$ with
intensity
$$
\lambda_t=\left\{
\begin{array}{ll}
\alpha t^{\alpha -1} &\mbox{ if }\; t>0\\
\beta |t|^{\beta -1} &\mbox{ if }\;  t<0
\end{array}
\right. .
$$
\end{cor}
The power of the martingale method can be seen when one wants to
obtain the asymptotic joint distribution of several $N^{(n)}(q)$.
\begin{thm}\label{theo.qqq}
Let  $0\le q_1<q_2<\ldots<q_k\le1$ and assume that
(\ref{def.univ.shift.rv.measure}) holds for each $x_{q_i},
i=1,...,k$, with $\alpha_i$ and $\beta_i$, respectively. Then
$$
\langle N^{(n)}(q_1), N^{(n)}(q_2),\ldots,
N^{(n)}(q_k)\rangle\longrightarrow_{\cal D} \langle
N(1),\ldots,N(k)\rangle\; ,
$$
where $\langle N(1),\ldots,N(k)\rangle$ is a $k$-variate Poisson
process on $\bbR$ with independent components and marginal
intensities $\lambda(i)$, $i=1,\ldots,k$, given by
$$
\lambda_t(i)=\left\{
\begin{array}{ll}
\alpha t^{\alpha_i -1} &\mbox{ if }\; t>0\\
\beta |t|^{\beta_i -1} &\mbox{ if }\;  t<0
\end{array}
\right. .
$$
\end{thm}
{\bf Proof:} For clarity, we will consider only the case $k=2$ and
verify the conditions of Theorem \ref{multichannel}. The general
result follows in a straightforward manner.

We begin by observing that it suffices to show joint convergence of
$$\langle N^{(n)}(q_1), N^{(n)}(q_2)\rangle$$ for all $t\in [-K,K]$
for any arbitrary finite constant $K$. Assume that $F$ is regularly
varying on the  right and left of $x_{q_i}$ with index $\alpha_i$
and $\beta_i$, respectively, $i=1,2$. As before, define $a_n^{(i)}$
 and $b_n^{(i)}$, $i=1,2$ according to (\ref{shifted.threshold}).
 Choose $M$ large enough that for $n\geq M$,
$[x_{q_1}-Kb_n^{(1)},x_{q_1}+Ka_n^{(1)}]$  and
$[x_{q_2}-Kb_n^{(2)},x_{q_2}+Ka_n^{(2)}]$ do not intersect. This
will ensure that those points $Y_j$ which are jump points of
$N^{(n)}(q_1)$ are not jump points of $ N^{(n)}(q_2)$ and vice
versa.

Consider for $i=1,2$, $1\le j\le n$ the single jump point process
$$J^{(n,j)}(q_i)=\delta_{(a_n^{(i)})^{-1}[Y_j-x_{q_i}]}I[Y_j\ge x_{q_i}]
+\delta_{(b_n^{(i)})^{-1}[Y_j-x_{q_i}]}I[Y_j< x_{q_i}].$$ It is
adapted to $\F=\{\F_{A_t}:-K\leq t\leq K\}$, where
$$\F_{A_t}=
\left\{
\begin{array}{ll}
\sigma\{\bbI_{\{Y_j\in
 [x_{q_i},x_{q_i}+ta_n^{(i)}] \}}, i=1,2; j=1,...,n\}&
\mbox{if } t\geq0\\
\sigma\{\bbI_{\{Y_j\in
 [x_{q_i}+tb_n^{(i)},x_{q_i  }]\}}, i=1,2; j=1,...,n\}&
\mbox{if } t\leq 0.
\end{array} \right.
$$
We will compute a *-compensator $\overline{J}^{(n,j)}(q_i)$ of the
single jump process $J^{(n,j)}(q_i)$. We consider only $0<t<K$ as
the argument for $t<0$ is similar. Let
$$U_t= [x_{q_1}-Kb_n^{(1)},x_{q_1}+ta_n^{(1)}]\cup
  [x_{q_2}-Kb_n^{(2)},x_{q_2}+ta_n^{(2)}].$$
Then for $C=(t,t']\in\C$, it follows that $\bbI_{\{Y_j\in
U_t\}}\in\G_C^*$ and so heuristically, the compensator
$\overline{J}^{(n,j)}(q_i)$ satisfies
$$\overline{J}_{dt}^{(n,j)}(q_i)=\frac{\bbI_{\{Y_j \in
U_t^c\}}dF(x_{q_i}+a_n^{(i)}t)}{1-F(U_t)},$$ provided that
$[x_{q_1}-Kb_n^{(1)},x_{q_1}+Ka_n^{(1)}]$  and
$[x_{q_2}-Kb_n^{(2)},x_{q_2}+Ka_n^{(2)}]$ are disjoint intervals.
 Using arguments similar to
those in Ivanoff and Merzbach (2000) it is   straightforward to
verify that for $n\geq M$ the
*-compensator $\Lambda^{(n)}(i)$ of $N^{(n)}(q_i)$ is
\begin{equation}\label{eq3}
\Lambda^{(n)}_{t}(i)=\sum_{j=1}^n\int_0^t\frac{\bbI_{\{Y_j \in
U_s^c\}}dF(x_{q_i}+a_n^{(i)}s)}{1-F(U_s)}.
 \end{equation}
Exactly as in the comments leading to Corollary \ref{theo.q} we have
$
 \bbE[\Lambda^{(n)}_{t}(i)]\sim t^{\alpha_i}
$ for the appropriate constant $\alpha_i$, since $F$ is slowly
varying on the right at $x_{q_i}$.

The argument that $\bbE[(\Lambda^{(n)}_{t}(i))^2]\rightarrow
(t^{\alpha_i})^2$ is also similar to that used in the proof of
Theorem \ref{theo.1dim}. Also, $N^{(n)}(q_i)$ is square integrable
with bounded  second moments (uniformly in $n$). Thus the conditions
of Theorem \ref{multichannel} have been satisfied and the result
follows.

\subsection{Multivariate case}\label{sec.mult}
Let $\{ Y_n\}_{n\ge 1}$ be a sequence of i.i.d. ${\matR}^d$-valued
random variables with continuous distribution $F$. Following the
pattern of the previous section, we may obtain a point process limit
if the regular variation index at ${\bf u}$ for $F$ depends on the
choice of orthant. To be precise, let $O_k$ be the $k$th orthant and
${\bf e}_k$ its associated unit vector, $k=1,\ldots,2^d$. Then $F$
is regularly varying at ${\bf u}$ from orthant ${\bf u}+O_k$, with
index $\alpha_k$ and rate $W_k$ if for ${\bf t}\in O_k$
\begin{equation}\label{eq.mev.measure.mult.shifted}
\lim_{ x\searrow 0}\frac{P_F(({\bf u},x{\bf t}+{\bf u}])}{P_F(({\bf
u},x{\bf e}_k+{\bf u}])}=W_k({\bf t})\; .
\end{equation}
The function $W_k$ is homogeneous of order $\alpha_k$, i.e. $W(s{\bf
t})=s^{\alpha_k}W({\bf t})$, see e.g. Resnick (1987).

We define $N^{(n)}$ in analogy to (\ref{eq2}), i.e.
$$
N^{(n)}=\sum_{k=1}^{2^d}\sum_{i=1}^n\delta_{a_{k,n}^{-1}[{ Y}_i-{\bf
u}]}I[{Y}_i\in O_k'],
$$
where $O'_{k}=O_k+{\bf u}$. More generally, if ${\bf u}_j\in
{\matR}^d$, $j=1,\ldots,m$, then we may define
\begin{equation}\label{eqn.multi}
N^{(n)}(j):=N^{(n)}({\bf
u}_j)=\sum_{k=1}^{2^d}\sum_{i=1}^n\delta_{a_{k,n}^{-1}[{ Y}_i-{\bf
u}_j]}I[{Y}_i\in O_{k,j}']\; ,
\end{equation}
where $O'_{k,j}=O_k+{\bf u}_j$.
\begin{thm}\label{theo.mult.multiv}
Assume that the orthant-wise regular variation conditions
(\ref{eq.mev.measure.mult.shifted}) are satisfied at ${\bf u}_{j},
j=1,...,m$, ${\bf x}_j\in {\matR}^d$. For each $j$, let $
N^{(n)}(j)$ denote the ${\matR}^d$-indexed point process of
(\ref{eqn.multi}). Then
$$
\langle N^{(n)}(1), N^{(n)}(2),\ldots, N^{(n)}(m)\rangle
\longrightarrow_{\cal D} \langle N(1), N(2),\ldots, N(m)\rangle
$$
where $\langle N(1), N(2),\ldots, N(m)\rangle$ is a vector of
independent Poisson processes, where the $j$th component process is
parameterized by ${\matR}^d$ and its mean measure is given
orthant-wise by the regular variation rates of $F$ at the
corresponding ${\bf u}_j$.
\end{thm}

Examples of regularly varying distributions are readily constructed.
One source of examples are distributions based on copulas as
described, for example, in Nelsen (1999) and Section \ref{mex}
below.


\section{Applications}\label{sec.applic}

\begin{remark}\label{rem.f}{\rm
The values of $a_n$ depend on the exact asymptotic behaviour of
the density at $x$, and certainly will not be known in general. In
our applications we consider only the special cases where the
slowly varying function $\ell(n)$ is in fact a constant, although
unknown. We then apply Theorem~\ref{theo.1dim} with the scaling
values equal to $n^{-1/\alpha}$.  We can define a compensator,
$\widetilde{\Lambda}^{(n)}$,  by (\ref{comp.empirical}) where
$a_n$ is replaced by $n^{-1/\alpha}$, and the relation to the
original definition is given by
$$
\widetilde{\Lambda}^{(n)}_t =
\Lambda^{(n)}_{\left({n^{-1/\alpha}}/{a_n}\right)t}.
$$
Since $\lim_{n\rightarrow
\infty}\left({n^{-1/\alpha}}/{a_n}\right)=\omega\in{\matR}$, then
$$
\widetilde{\Lambda}^{(n)}_t =
\Lambda^{(n)}_{\left({n^{-1/\alpha}}/{a_n}\right)t}\rightarrow
\Lambda_{\omega t}
$$
and we have convergence of the empirical point process to a Poisson
process with intensity $\omega\alpha t^{\alpha-1}$.  For example, if
the density of $Y$ at 0 is $8$, then the weak limit of
$N^{(n)}=\sum_{i=1}^n \delta_{n^{-1/\alpha}Y_i}$ is a Poisson
process with intensity $8\alpha t^{\alpha-1}$.  The corresponding
changes to the other theorems of the previous section are
immediate.}
\end{remark}

\subsection{Local Density Estimation}
Consider a sample $\{ Y_1, Y_2, \ldots, Y_n\}$ with common
marginal differentiable distribution $F$ on $[0,1]$, and assume
that its density $f$ is positive on the range $[0,1]$.  Let $F_n$
denote the empirical distribution and define $[t]_n^+=Y_{(k+1)}$
and $[t]_n^-=Y_{(k)}$ by $Y_{(k)}\le t<Y_{(k+1)}$. We put
$[t]_n^+=0$ if $[t]_n^+<Y_{(1)}$ and $[t]_n^-=0$ if
$[t]_n^+>Y_{(n)}$. A naive nearest-neighbour estimator of the
density at t is given by
\begin{equation}\label{nearest}
\widehat{f}(n,t)=\frac{1}{n}\left/\left(([t]_n^+
-t)+(t-[t]_n^-)\right)\right. .
\end{equation}
Additional information on nearest-neighbour density estimates can be
found in H\"{a}rdle (1990) or Silverman (1992), including comments
on performance, and modifications.

For $t\in(0,1)$ the fact that $F$ is differentiable and that $f$ is
positive (i.e. $F$ is of regular variation index $\alpha=1$ at $t$)
allows us to write
\begin{equation}
\widehat{f}(n,t)/f(t)= 1\left/f(t)\left(n([t]_n^+
-t)+n(t-[t]_n^-)\right)\right. \stackrel{\mathcal D}{\rightarrow}
1/f(t)(E_1+E_2) \label{eqn.naive}
\end{equation}
where $E_1$ and $E_2$ are independent exponential variables of mean
$1/f(t)$.  This convergence follows from Corollary~\ref{theo.q} and
the continuous mapping theorem, and follows the pattern set for
extreme value processes as given in Resnick (1987). As each limiting
Poisson process has a constant rate function equal to $f(t)$, the
distance from $t$ to the first point has an exponential distribution
with mean $1/f(t)$.  Since such an exponential variable can be
written as the product of $1/f(t)$ and an exponential of mean 1, and
the sum of two independent mean 1 exponentials is a $\Gamma(2,1)$
variable, we have identified the limiting distribution of
$\widehat{f}(n,t)/f(t)$ as Inverse Gamma, $\Gamma^{(-1)}(2,1)$. The
mode, mean and variance of an Inverse Gamma density of parameters
$(\alpha,\beta)$ are $\beta/(\alpha+1)$, $\beta/(\alpha-1)$ (for
$\alpha>1$) and $\beta^2/((\alpha-1)(\alpha-2))$ (for $\alpha>2$),
respectively. Thus we see that this naive estimator of $f(t)$ has
mode ${f(t)}/{3}$, mean $f(t)$ and infinite variance.

This development can be easily extended to estimators based on the
$k$ lower nearest neighbours and $k$ upper nearest neighbours.  As
above, asymptotically the spacings between consecutive neighbours
are independent exponential variables with mean $1/f(t)$.  The
asymptotic joint density is the product of $2k$ exponentials, and
the sufficient statistic is just the total distance from the lower
$k$th-nearest neighbour of $t$, $[t]_{n}^{-k}$, to the upper
$k$th-nearest neighbour, $[t]_{n}^{+k}$.

\begin{cor} The
asymptotically uniformly minimum variance unbiased estimator ($k>1$)
is
$$
\widehat{f}_k(n,t)=\frac{(2k-1)/n}{[t]_{n}^{+k}-[t]_{n}^{-k}},
$$
and $\widehat{f}_k(n,t)/f(t)$ has an asymptotic
$\Gamma^{(-1)}(2k,1)$ density.
\end{cor}

Using this result we can consequently compute approximate confidence
intervals for $f(t)$ or construct tests. If $k$ is fixed,
Theorem~\ref{theo.qqq} also identifies the limiting distribution of
$$\langle \widehat{f}_k(n,t_1), \widehat{f}_k(n,t_2), \ldots,
\widehat{f}_k(n,t_m) \rangle$$ as given by a vector of $m$
independent scaled inverse Gamma variables.  Consequently we can
obtain the limiting distribution of expressions such as approximate
integrals,
$$
\widehat{E}(g(Y))=\sum_{i=1}^m g(t_i)(\widehat{f}_n(t_i)),
$$
even for arbitrary dimension (Theorem~\ref{theo.mult.multiv}) with
appropriate norming.

\begin{remark}\rm
On the other hand, we see that $\widehat{f}_k(n,t)/f(t)$ still has
an Inverse Gamma distribution, but with finite variance for $k\ge
1$. It has asymptotic variance $1+1/(2k-2)$, and so remains
inherently random regardless of the fixed number of nearest
neighbours used in the estimate. Nearest-neighbour methods have
become popular in data mining, classification and computing
applications, and rapid algorithms exist for finding the $k$
nearest neighbours to a point $t$ even in high dimensions. The
above discussion shows that even in highly regular cases, the best
$k$-nearest-neighbour density estimate will not converge in
probability to the desired limit, and remains random.
\end{remark}

As an example of a test that can be constructed using the results of
this paper, we consider the null hypothesis that $F$ is regularly
varying as $\omega t$ from the right at $0$ for some $\omega>0$
(e.g. $F^\prime(0)=\omega>0$). We take the alternative to be where
$F$ varies as $\zeta t^2$ from the right for some $\zeta>0$ (i.e.
$F^\prime(0)=0$). The maximum likelihood under the null hypothesis
is proportional to $([t]_n^{+k})^{-k}$, and that under the
alternative  proportional to
$([t]_n^{+k})^{-k}\times\prod_{i=1}^k([t]_n^{+i}/[t]_n^{+k})$.

\begin{cor}
 The likelihood ratio test based on the $k\ge 2$ upper nearest neighbours
rejects when $$\prod_{i=1}^{(k-1)}([t]_n^{+i}/[t]_n^{+k})$$ is too
large. Under the null hypothesis, the distribution of this product
is given by the product of $k-1$ independent uniform variables on
$[0,1]$. \end{cor}

When $k=2$ we obtain an intuitively reasonable test that rejects
when the distance from $0$ to $[t]_n^{+1}$ is much larger than that
from $[t]_n^{+1}$ to $[t]_n^{+2}$, and so indicates the presence of
a ``gap'' in the distribution.

\subsection{Multivariate extremes}\label{mex} Let $\{(Y_{n1},Y_{n2})\}_{n\ge
1}$ be an i.i.d. sequence of bivariate random vectors. To focus on
the bivariate dependence structure rather than the marginal
distributions, we assume that $(Y_{11},Y_{12})$ has a copula $C$ and
standard uniform marginals, see Nelsen (1999). We want to
characterize
$$
{\matP}(Y_{11}>1-xt_1,Y_{12}>1-xt_2)
$$
as $x\searrow 0$. If $Y_{11}$ and $Y_{12}$ are independent, then the
above probability factors and we can apply standard extreme value
methods (e.g. Resnick, 1987) to the marginals.  However, if $Y_{11}$
and $Y_{12}$ are dependent but the maxima are asymptotically
independent then the extreme value methods fail; see Fougeres (2004)
for a general discussion of this problem. For most known families of
copulas which have the asymptotic independence property, we have
(cf. Hefferman, 2000)
\begin{equation}
{\matP}(Y_{11}>1-xt_1,Y_{12}>1-xt_2)\sim c x^2.\label{eqn:cop}
\end{equation}
By the results of this paper the appropriate scaling to obtain a
point process limit for the joint extremes is $a_n=n^{-1/2}$, and
not the $a_n=n^{-1}$ that would be used to normalize the marginal
variables individually. Note, moreover, that the methods of this
paper are ``dimension free'', and so we can address multivariate
copulas of any dimension.

Further we can address the joint extreme value behaviour of
copulas with the asymptotic independence property but where
(\ref{eqn:cop}) is not satisfied. Consider the case when $C$ is
the bivariate normal copula with correlation $\rho\in (0,1]$ --
i.e. $C(x,y)$ is given by a joint normal distribution function at
$(\Phi^{-1}(x),\Phi^{-1}(y))$ with standard marginals and
correlation $\rho$. We have
$$
{\matP}(Y_{11}>1-xt_1,Y_{12}>1-xt_2)\sim x^{2/(1+\rho)}g(t_1,t_2)
$$
for a function $g$ as $x\searrow 0$, and so $$ \lim_{x\searrow
0}\frac{{\matP}(Y_{11}>1-xt_1,Y_{12}>1-xt_2)}
{{\matP}(Y_{11}>1-x,Y_{12}>1-x)}=\frac{g(t_1,t_2)}{g(1,1)}\; ,
$$
for $t_1,t_2\geq 0$. For ${\bf u}=(1,1)$ and $t_1,t_2\leq 0$,
formula (\ref{eq.mev.measure.mult.shifted}) is satisfied with
$$W(t_1,t_2)=\frac{g(-t_1,-t_2)}{g(1,1)}.$$ Applying the results of
Section \ref{sec.mult}, we can characterize the asymptotics of joint
extremes for a normal copula.
\begin{cor}\label{corcopula}
Assume that $\{{\bf Y}_n=(Y_{n1},Y_{n2})\}_{n\ge 1}$ are
independent, have a common normal copula of parameter $\rho$ and
uniform marginals. Then for ${\bf u}=(1,1)$ and
$a_n=n^{-(1+\rho)/2}$,
$$N^{(n)}=\sum_{i=1}^n\delta_{a_n^{-1}[{\bf Y}_i-{\bf u}]}$$
converges to a Poisson process on ${\matR}_{-}^2$ with mean measure
$W(\cdot,\cdot)$.
\end{cor}

\ack The research of the first two authors was supported by grants
from the Natural Sciences and Engineering Research Council of
Canada. This paper was written while Rafa{\l} Kulik was a
postdoctoral fellow at the University of Ottawa.


%
%
%
%

\end{document}